\documentclass[11pt]{article}

\usepackage[algoruled,linesnumbered]{algorithm2e}
\usepackage{algorithmic}
\newcommand{\aut}{{\rm Aut}}
\usepackage{amsmath,amssymb}
\usepackage{graphicx}
\usepackage{epstopdf}
\usepackage{caption}
\usepackage{color}
\usepackage{dcolumn}
\usepackage{bm}
\usepackage[numbers,super,comma,sort&compress]{natbib}

\definecolor{background-color}{gray}{0.98}

\title{Predicting Melting Points by the Graovac-Pisanski Index}
\author{Matev\v z \v Crepnjak \thanks{University of Maribor, Faculty of Natural Sciences and Mathematics, Koro\v ska cesta 160, 2000 Maribor, Slovenia; University of Maribor, Faculty of Chemistry and Chemical Engineering, Smetanova ulica 17, 2000 Maribor, Slovenia; University of Primorska, Andrej Maru\v si\v c Institute, Muzejski trg 2, 6000 Koper, Slovenia}, Niko Tratnik\thanks{University of Maribor, Faculty of Natural Sciences and Mathematics, Koro\v ska cesta 160, 2000 Maribor, Slovenia}, Petra \v Zigert Pleter\v sek \thanks{University of Maribor, Faculty of Natural Sciences and Mathematics, Koro\v ska cesta 160, 2000 Maribor, Slovenia; University of Maribor, Faculty of Chemistry and Chemical Engineering, Smetanova ulica 17, 2000 Maribor, Slovenia}}

\begin{document}

\maketitle

\begin{abstract}
Theoretical molecular descriptors alias  topological indices are a convenient means for expressing in a numerical form the chemical structure encoded in a molecular graph. The structure descriptors derived from molecular graphs are widely used in Quantitative Structure-Property Relationship (QSPR) and  Quantitative Structure-Activity Relationship (QSAR). 

In this paper, we are interested in  the Graovac-Pisanski index (also called modified Wiener index)  introduced in 1991 by Graovac and Pisanski, which encounters beside the distances in a molecular graph also its symmetries. In the QSPR analysis we first calculate the  Graovac-Pisanski index for different families of hydrocarbon molecules  using a simple program  and then we show  a correlation with the melting points of considered molecules. We show that the melting points of the alkane series can be very effectively predicted by the Graovac-Pisanski index and  for the rest of considered molecules (PAH's and octane isomers) the regression models are different, but we establish some  correlation with the melting points for them as well.

\end{abstract}

\clearpage


  \makeatletter
  \renewcommand\@biblabel[1]{#1.}
  \makeatother

\bibliographystyle{apsrev}

\renewcommand{\baselinestretch}{1.5}
\normalsize

\clearpage

\section*{\sffamily \Large INTRODUCTION} 

Theoretical molecular descriptors (also called topological indices) are graph invariants that play an important role in chemistry, pharmaceutical sciences, materials science and engineering, etc. The value of a molecular descriptor must be independent of the particular characteristics of the molecular representation, such as atom numbering or labeling. We model molecules of hydrocarbons by the corresponding molecular graph, where the vertices are the carbon atoms and the edges of the graph are the bonds between them.

One of the most investigated topological indices is the Wiener index introduced in 1947 \cite{wiener}. This index is defined as the sum of distances between all the pairs of vertices in a molecular graph. Wiener showed that the Wiener index is closely correlated with the boiling points of alkane molecules. 
In order to take into account also the symmetries of a molecule, Graovac and Pisanski in 1991 introduced the modified Wiener index \cite{graovac}. However, the name modified Wiener index was later used for different variations of the Wiener index and therefore, Ghorbani and Klav\v zar suggested the name Graovac-Pisanski index\cite{ghorbani}, which is also used in this paper.

Very recently, the Graovac-Pisanski index of some molecular graphs and nanostructures was extensively studied\cite{ashrafi_koo_diu,ashrafi_sha,
koo_ashrafi3,koo_ashrafi,koo_ashrafi2,sha_ashrafi}.  Moreover, the closed formulas for the Graovac-Pisanski index of zig-zag nanotubes were computed\cite{tratnik}.  The only known connection  of the Graovac-Pisanski index with some molecular properties is the  correlation with the topological efficiency \cite{ashrafi_koo_diu1}. Therefore, it was pointed out by Ghorbani and Klav\v zar \cite{ghorbani}  that the  QSPR or QSAR analysis should be performed in order to establish correlation with some other physical or chemical properties of molecules.

On the other hand, there is no known molecular descriptor well correlated to the melting points of molecules, since graph-theoretical abstraction in most cases disregards many information that are relevant for the value of the melting point\cite{ro-king,vuk-gas}.

The symmetries of a molecule play an important role in the process of melting\cite{pinal}. The more symmetrical the molecules are, easier it is for them to stack together. Consequently, the fewer spaces there are between them and so  the melting  point is expected to be  higher. Therefore, since the Graovac-Pisanski index considers the symmetries of a molecule, in this paper we investigate its correlation with the melting points of molecules.

We work with the data set of molecules proposed by the International Academy of Mathematical Chemistry, in particular alkane series, polyaromatic hydrocarbons (PAH's) and octane isomers. We show that for alkane series the melting point is very well correlated with the Graovac-Pisanski index and for PAH's the correlation is a little bit weaker. For octane isomers, the melting points are good correlated with the number of symmetries. However, we conclude that for alkanes the size of a molecule contributes the most to its melting point.

\section*{\sffamily \Large THE GRAOVAC-PISANSKI INDEX}

A \textit{graph} $G$ is an ordered pair $G = (V, E)$ of a set $V$ of \textit{vertices} (also called nodes or points) together with a set $E$ of \textit{edges}, which are $2$-element subsets of $V$ (more information about some basic concepts in graph theory can be found in a book written by West\cite{west}). Having a molecule, if we represent atoms by vertices and bonds by edges, we obtain a \textit{molecular graph}.

The graphs considered in this paper are all finite and connected. The {\em distance} $d_G(x,y)$ between vertices $x$ and $y$ of a graph $G$ is the length of a shortest path between vertices $x$ and $y$ in $G$ (we often use $d(x,y)$ for $d_G(x,y)$). 
\bigskip

\noindent
The {\em Wiener index} of a graph $G$ is defined as $\displaystyle{W(G) = \frac{1}{2} \sum_{u \in V(G)} \sum_{v \in V(G)} d_G(u,v)}$. Moreover, if $S \subseteq V(G)$, then $\displaystyle{W(S) = \frac{1}{2} \sum_{u \in S} \sum_{v \in S} d_G(u,v)}$.
\bigskip

\noindent  
An \textit{isomorphism of graphs} $G$ and $H$ with $|E(G)|=|E(H)|$ is a bijection $f$ between the vertex sets of $G$ and $H$, $f: V(G)\to V(H)$,
such that for any two vertices $u$ and $v$ of $G$ it holds that if $u$ and $v$ are adjacent in $G$ then $f(u)$ and $f(v)$ are adjacent in $H$. When $G$ and $H$ are the same graph, the function $f$ is called an \textit{automorphism} of $G$. The composition of two automorphisms is again an automorphism, and the set of automorphisms of a given graph $G$, under the composition operation, forms a group ${\aut}(G)$, which is called the \textit{automorphism group} of the graph $G$.
\bigskip

\noindent 
The \textit{Graovac-Pisanski index} of a graph $G$, $GP(G)$, is defined as
$$GP(G) = \frac{|V(G)|}{2 |{\aut}(G)|} \sum_{u \in V(G)} \sum_{\alpha \in {\aut}(G)} d_G(u, \alpha(u)).$$

\noindent 
Next, we mention some important concepts of group theory. If $G$ is a group and $X$ is a set, then a \textit{group action} $\phi$ of $G$ on $X$ is a function $\phi :G \times X \to X$
that satisfies the following: $\phi(e,x) = x$ for any $x \in X$ (where $e$ is the neutral element of $G$) and $\phi(gh,x)=\phi(g,\phi(h,x))$ for all $g,h \in G$ and $x \in X$. The \textit{orbit} of an element $x$ in $X$ is the set of elements in $X$ to which $x$ can be moved by the elements of $G$, i.e.\ the set $\lbrace \phi(g,x) \, | \, g \in G \rbrace$. If $G$ is a graph and ${\aut}(G)$ the automorphism group, then $\phi: {\aut}(G) \times V(G) \to V(G)$, defined by $\phi(\alpha,u) = \alpha(u)$ for any $\alpha \in {\aut}(G)$, $u \in V(G)$, is called the \textit{natural action} of the group ${\aut}(G)$ on $V(G)$.

\noindent
It was shown by Graovac and Pisanski\cite{graovac} that if $V_1, \ldots, V_t$ are the orbits under the natural action of the group ${\aut}(G)$ on $V(G)$, then
\begin{equation*}
\label{formula}
GP(G) = |V(G)| \sum_{i=1}^t \frac{1}{|V_i|}W(V_i).
\end{equation*}

\begin{figure}[h!]
\centering
\includegraphics[width=0.5\columnwidth,keepaspectratio=true]{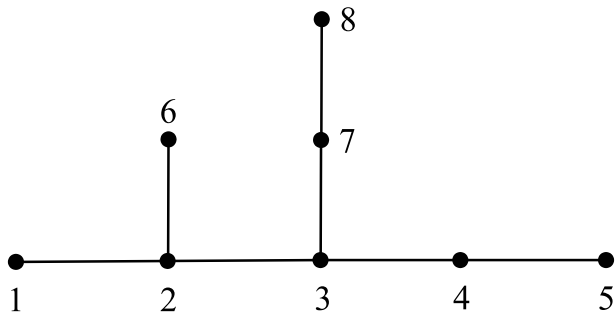}
\caption{\label{primer_drevo} Molecular graph $G$ of 2-methyl-3-ethyl-pentane.}
\end{figure}

As an example, we calculate the Graovac-Pisanski index for one of the octane isomers, more precisely for 2-methyl-3-ethyl-pentane, see Figure \ref{primer_drevo}. 

We do the calculation in two different ways. First, we calculate directly by the definition.
There are all together 3 nontrivial (different from the identity) automorphisms  of the considered molecular graph $G$:
$$\alpha_1 = (1\ 6)(2)(3)(4)(5)(7)(8),\ \alpha_2 = (1)(2)(3)(4\ 7)(5\ 8)(6),\ \alpha_3 = (1\ 6)(2)(3)(4\ 7)(5\ 8).$$
Therefore,
\begin{align*} 
GP(G) & =  d(1,\alpha_1(1)) + d(6,\alpha_1(6)) + d(4,\alpha_2(4)) + d(7,\alpha_2(7)) + d(5,\alpha_2(5)) + d(8,\alpha_2(8))  \\
			& +  d(1,\alpha_3(1)) + d(6,\alpha_3(6)) + d(4,\alpha_3(4)) + d(7,\alpha_3(7)) + d(5,\alpha_3(5)) + d(8,\alpha_3(8)) \\
			& = 2+2+2+2+4+4+2+2+2+2+4+4 \\
			& = 32. 
\end{align*}

Now, we calculate the index by using orbits.
The vertex set is partitioned into 5 orbits
$$V_1=\{1,6\},\,V_2=\{2\},\,V_3=\{3\},\,V_4=\{4,7\},\,V_5=\{5,8\}$$
and the Graovac-Pisanski index of 2-methyl-3-ethyl-pentane is calculated as
$$ GP(G)  =  |V(G)| \sum_{i=1}^5 \frac{1}{|V_i|}W(V_i)
      =8\left( \frac{2}{2}+\frac{2}{2}+\frac{4}{2}\right)=32\,. $$

\section*{\sffamily \Large COMPUTATIONAL DETAILS}


In this section we present an algorithm which was used to compute the Graovac-Pisanski index and the number of automorphisms of a graph. The algorithm contains two special functions, i.e. {\tt calculateAutomorphisms} and {\tt calculateDistances}.

Let $G$ be a graph represented by a adjacency matrix with vertices $1,2,\ldots,n$. The function  {\tt calculateAutomorphisms} determines all the automorphisms of graph $G$ and saves them in a set $A$. One possibility is to go through all the permutations of the set $\lbrace 1,2,\ldots, n \rbrace$ and check  whether it is an automorphism. Note that for the graph automorphism problem (which is the problem of testing whether a graph has a nontrivial automorphism)  it is still unknown whether it has a polynomial time algorithm or it is NP-complete\cite{lubiw}.

The function  {\tt calculateDistances} computes the matrix $M$ from the adjacency matrix of $G$. The element $M_{i,j}$, $i,j \in \lbrace 1,2,\ldots, n \rbrace$ of $M$ represents the distance between vertices $i$ and $j$ in $G$. Note that this algorithm is known as Floyd-Warshall algorithm \cite{floyd} and has the time complexity $O(n^3)$.

\begin{algorithm}[H]\label{alg:edini}
\SetKwInOut{Input}{Input}\SetKwInOut{Output}{Output}
\DontPrintSemicolon

\Input{Graph $G$ with vertices $1,2,\ldots,n$.}
\Output{$GP(G)$ and $|\aut(G)|$}

 \SetKwFunction{CA}{calculateAutomorphisms}
 \SetKwFunction{CD}{calculateDistances}
 $A \leftarrow$ \CA($G$)\;
 $M \leftarrow$ \CD($G$)\;
 $X \leftarrow 0$\;
 
 \For{\bf{each} $\alpha \in A$}{
 \For{$i=1$ {\rm to} $n$}{
 	$X \leftarrow X + M_{i,\alpha(i)}$ \;
 	
 }}
 $GP(G) \leftarrow \frac{n}{2|A|}X$ \;
 $|\aut(G)| \leftarrow |A|$
 \caption{The Graovac-Pisanski index}
\end{algorithm}

\bigskip
We notice that the time complexity of the algorithm is not polynomial if we go through all the permutations. However, for many chemical graphs the set of all the automorphisms of a graph can be easily obtained by hand. In such a case, we can skip the first line of the algorithm and consequently, it becomes much more efficient.



\section*{\sffamily \Large RESULTS AND DISCUSSION}


\section*{\sffamily \Large Alkanes}
 
In this section, we have 31 alkane molecules, which are divided into two sets. The training set contains 26 molecules and the test set has 5 molecules (in Table \ref{tab1} the molecules in the test set are written in bold style).

\begin{center}
\begin{figure}[h!]
\centering
\includegraphics[width=0.9\columnwidth,keepaspectratio=true]{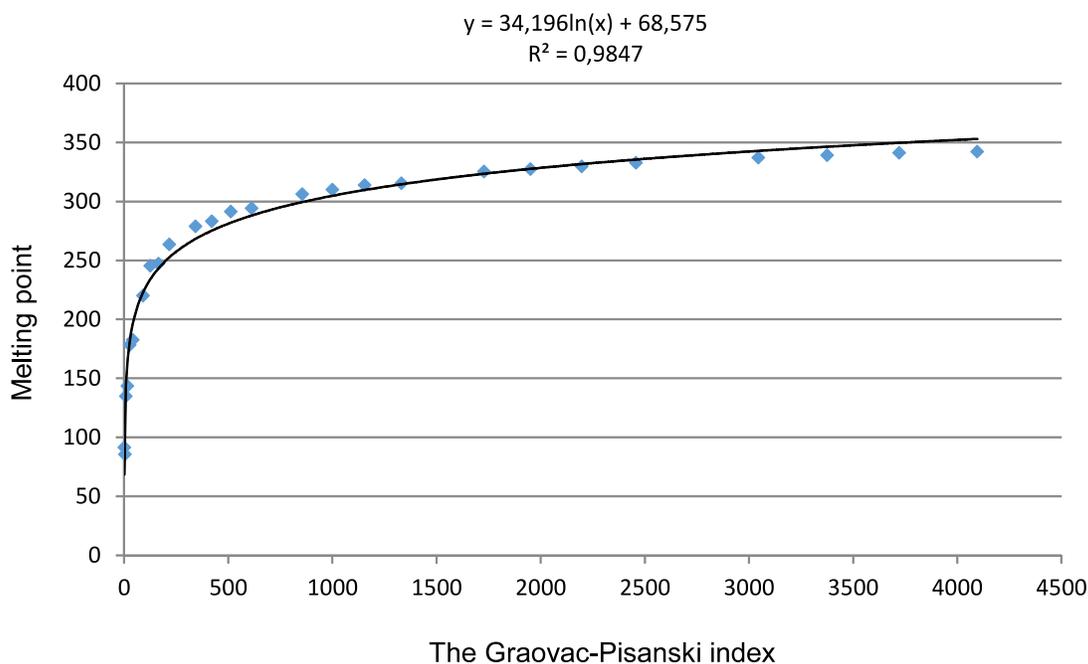}
\caption{\label{figu1} Nonlinear regression on the training set.}
\end{figure}

\end{center}

The data show that a logarithmic function $f(x)= a \ln x + b$ fits the best among all elementary functions. After performing nonlinear regression on the training set, we obtain $a=34,196$ and $ b=68,575$, see Figure \ref{figu1}. Therefore, the following equation can be used to predict the melting points of the alkane series
$$MP = 34,\!196 \ln GP + 68,\!575.$$

\noindent
We use the above formula to compute the predicted melting point $\widehat{MP}$ for every molecule from the test set, see Table \ref{tab2}. From Table \ref{tab2} we can see that the average error on the test set is less than 2\%. Moreover, the statistics shows very good correlation since $R^2=0,\!9847$.

\begin{center}
\begin{figure}[h!]
\centering
\includegraphics[width=0.9\columnwidth,keepaspectratio=true]{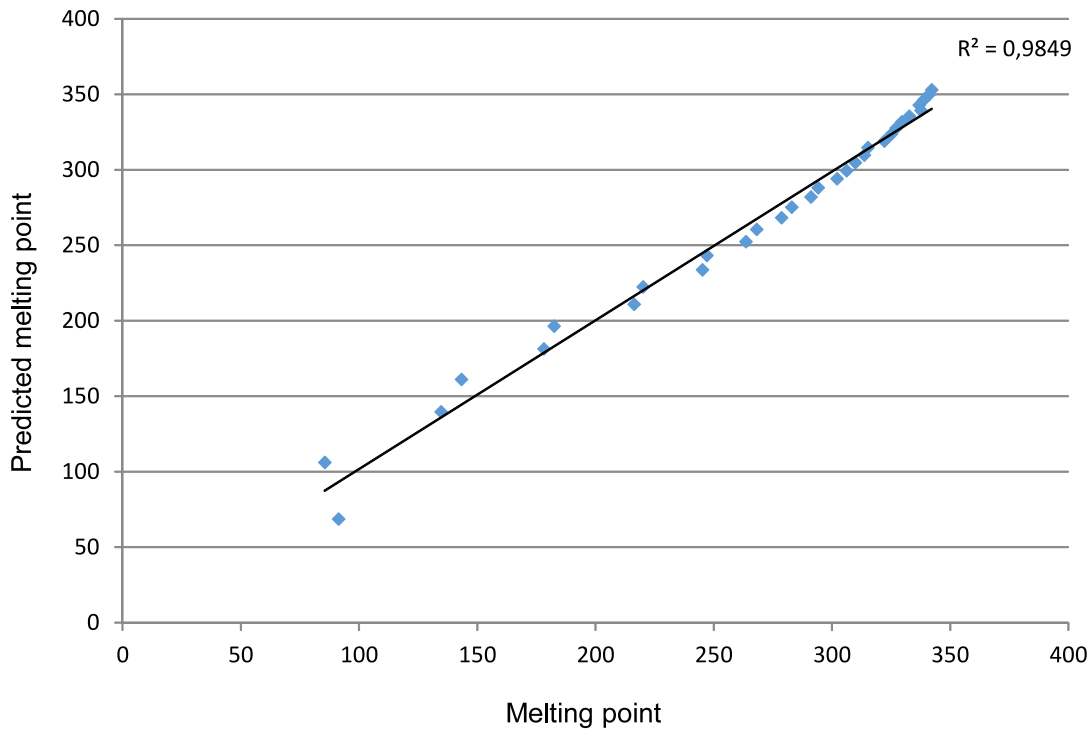}
\caption{\label{figu2} Melting points and predicted melting points for alkane molecules.}
\end{figure}

\end{center}

To conclude the section, we test the obtained formula on all the molecules from Table \ref{tab1} and we obtain Table \ref{tab3}. We can see that the average error is around 4\%. Figure \ref{figu2} shows the comparison between the melting points and predicted melting points. 

\section*{\sffamily \Large Polyaromatic hydrocarbons (PAH's)}

\begin{center}
\begin{figure}[h!]
\includegraphics[width=0.9\columnwidth,keepaspectratio=true]{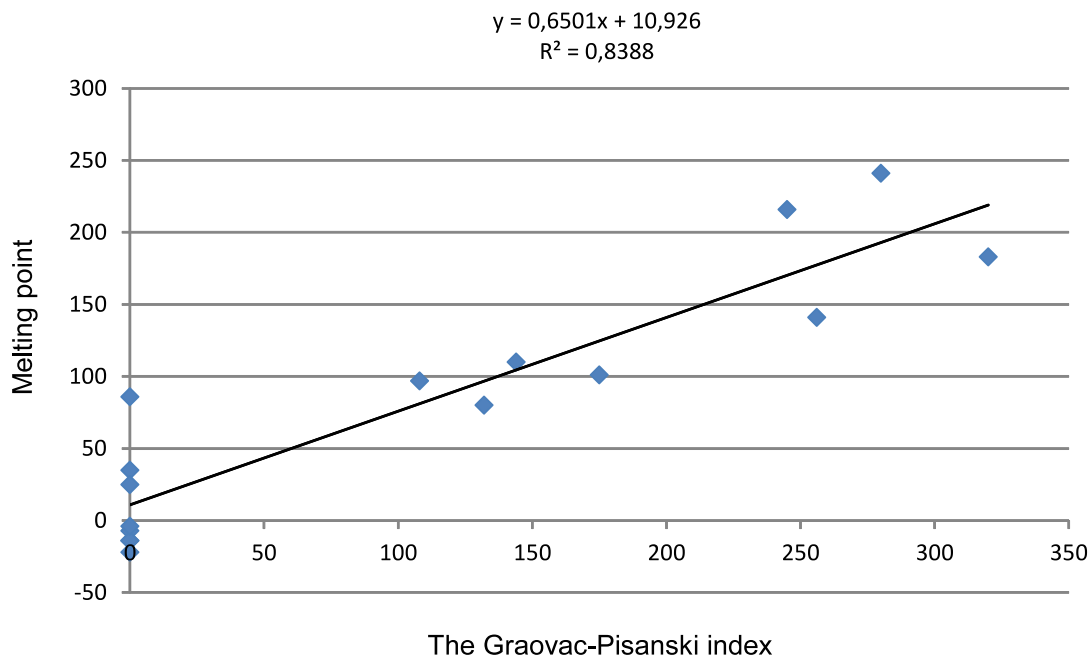}
\caption{\label{figu3} Linear regression on the training set for PAH's.}
\end{figure}
\end{center}

In this section, we consider 20 PAH molecules. In the first part, we divide these molecules into two sets (the training set contains 16 molecules and the test set has 4 molecules) and perform linear regression with respect to the Graovac-Pisanski index. In the second part, we improve the correlation by performing multilinear regression with respect to the Graovac-Pisanski index, the Wiener index and the number of automorphisms. The data for the PAH molecules is collected in Table \ref{tab4}. The linear regression results in the function (see Figure \ref{figu3})
$$MP = 0,\!6501 \, GP + 10,\!926.$$

\noindent
We use the above formula to compute the predicted melting point $\widehat{MP}$ for every molecule from the test set, see Table \ref{tab5}. From Table \ref{tab5} we can see that the average error on the test set is around 11\% and the statistics shows quite good correlation since $R^2=0,\!8388$.

To improve the correlation, we perform multilinear regression on the whole set of PAH's as described in the beginning of the section. It results in the formula

$$MP = -46,\!248 + 13,\!038\,(\# {\rm Aut}) + 0,\!446\,{GP} + 0,\!235\,{W}.$$

\noindent
The regression statistics (multiple $R$ is 0,946; 	$R^2$ is 0,894; adjusted $R^2$ is 0,874; and standard error is 30,665) shows better correlation than the linear regression.

\section*{\sffamily \Large 14 octane isomers}

As the last one, we consider the set of 14 octane isomers. However, all together there are 18 octane isomers, but for 4 of them, the data for the melting point was unavailable. We compute the Graovac-Pisanski index and the number of automorphisms for these molecules, see Table \ref{tab6}.

\begin{center}
\begin{figure}[h!]
\centering
\includegraphics[width=0.9\columnwidth,keepaspectratio=true]{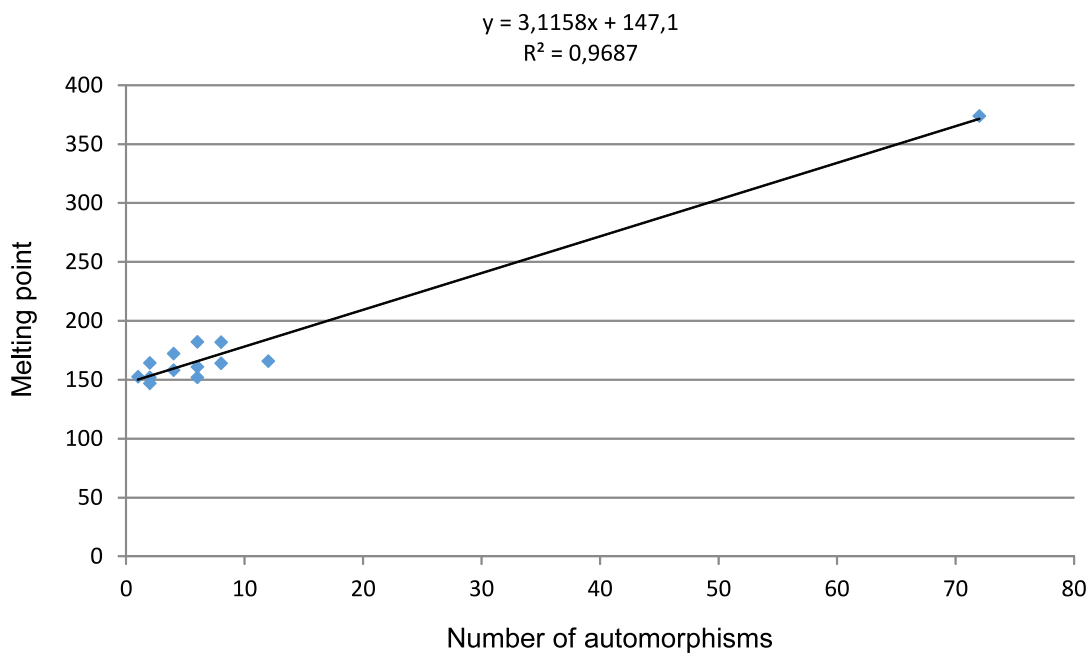}
\caption{\label{figu4} Linear regression between the number of automorphisms and the melting point for octane isomers.}
\end{figure}
\end{center}

It turns out that the correlation between the Graovac-Pisanski index and the melting point is not that good ($R^2$ is 0,2423), but there is good correlation ($R^2$ is 0,9687) between the number of automorphisms and the melting point, see Figure  \ref{figu4}. However, we can see that the data for 2,2,3,3-tetramethylbutane is standing out. If we exclude this molecule from our observation, the correlation between the number of automorphisms and the melting points becomes much weaker (close to zero), but the correlation between the Graovac-Pisanski index and the melting points becomes slightly stronger ($R^2$ is 0,4537). Therefore, we can conclude that the size of an alkane molecule contributes the most to its melting point.

%
%

%
%
%
%


\section*{\sffamily \Large CONCLUSIONS}


In recent years revived Graovac-Pisanki index was considered and calculated in different ways for some families of graphs, such as fullerene graphs or catacondensed benzenoid graphs, but almost no chemical application was known so far. 
In this paper we use the    QSPR analysis and show that this index is    well correlated with the melting points of the alkane series and polyaromatic hydrocarbon molecules. Beside that the connection between the number of graph automorphisms  of the octane isomers and theirs  melting point is established. 
In some sense this result in similar to the seminal paper  in the field of molecular descriptors \cite{wiener} by Wiener from 1947, where it was shown that   the boiling points of the alkane series can be predicted from the Wiener index.  Therefore, these results might contribute to further development in the area of obtaining   the Graovac-Pisanski index for different molecular graphs.

\subsection*{\sffamily \large ACKNOWLEDGMENTS}

\noindent The authors Matev\v z \v Crepnjak and Petra \v Zigert Pleter\v sek acknowledge the financial support from the Slovenian Research Agency, research core funding No. P1-0285 and No. P1-0297, respectively. The author Niko Tratnik was financially supported by the Slovenian Research Agency.

\clearpage



 \begin{table}
 \centering
\begin{tabular}{|c|c|c|}\hline  
Alkane & $GP$ & $MP \,(K)$  \\ \hline
ethane	& 1 &	91,39	 \\ 
propane	& 3	& 85,45	 \\
butane	& 8	& 134,75	 \\
pentane	& 15 & 143,35	 \\
hexane	& 27	& 178,15	 \\
heptane	& 42	& 182,54	 \\
\textbf{octane}	& 64	& 216,3	 \\
nonane	& 90	& 220,15	 \\
decane	& 125	& 245,25	 \\
undecane	& 165	& 247,15 \\
dodecane	& 216	& 263,55	 \\
\textbf{tridecane}	& 273	& 268,15	 \\
tetradecane	& 343	& 278,65	 \\
pentadecane	& 420	& 283,05	 \\
hexadecane	& 512	& 291,15	 \\
heptadecane	& 612	& 294,15	 \\
\textbf{octadecane}	& 729	& 302,15	 \\
nonadecane	& 855	& 306,15	 \\
icosane	& 1000	& 309,85	 \\
henicosane	& 1155	& 313,65	 \\
docosane	& 1331	& 315,15	 \\
\textbf{tricosane}	& 1518	& 322,15	 \\
tetracosane	& 1728	& 325,15	 \\
pentacosane	& 1950	& 327,15	 \\
hexacosane	& 2197	& 329,55	 \\
heptacosane	& 2457	& 332,65	 \\
\textbf{octacosane}	& 2744	& 337,65	 \\
nonacosane	& 3045	& 336,85	 \\
triacontane	& 3375	& 338,95	 \\
hentriacontane	& 3720	& 341,05   \\
dotriacontane	& 4096	& 342,15   \\
 \hline
 
\end{tabular}
\caption{\label{tab1} The Graovac-Pisanski index $GP$ and the melting point $MP$ for 31 alkane molecules.}
\end{table}

\begin{table}
\centering
\begin{tabular}{|c|c|c|c|c|c|}\hline  
Alkane & $GP$ & $MP\,(K)$ & $\widehat{MP}$	& Residual & \% Residual \\ \hline

octane	& 64	& 216,3	& 210,792	& 5,508	& 2,546 \\

tridecane	& 273	& 268,15	& 260,396	& 7,754	& 2,891 \\

octadecane	& 729	& 302,15	& 293,984	& 8,166	& 2,703 \\

tricosane	& 1518	& 322,15	& 319,066	& 3,084	& 0,957 \\

octacosane	& 2744	& 337,65	& 339,311	& -1,661	& 0,492 \\
\hline
average				& & & & &	1,918 
\\ \hline  

\end{tabular}
\caption{\label{tab2} Data for 5 alkane molecules in the test set.}
\end{table}

\begin{table}
\centering
\begin{tabular}{|c|c|c|c|c|c|}\hline  
Alkane & $GP$ & $MP\,(K)$ & $\widehat{MP}$	& Residual & \% Residual \\ \hline
ethane	& 1 &	91,39	& 68,575 &	22,815 &	24,964 \\ 
propane	& 3	& 85,45	& 106,143	& -20,693	& 24,217 \\
butane	& 8	& 134,75	& 139,684	& -4,935	& 3,661 \\
pentane	& 15 & 143,35	& 161,179	& -17,829	& 12,438 \\
hexane	& 27	& 178,15	& 181,279	& -3,129	& 1,757 \\
heptane	& 42	& 182,54	& 196,388	& -13,848	& 7,586 \\
octane	& 64	& 216,3	& 210,792	& 5,508	& 2,546 \\
nonane	& 90	& 220,15	& 222,450	& -2,300	& 1,045 \\
decane	& 125	& 245,25	& 233,684	& 11,566	& 4,716 \\
undecane	& 165	& 247,15	& 243,178	& 3,972	& 1,607 \\
dodecane	& 216	& 263,55	& 252,388	& 11,162	& 4,235 \\
tridecane	& 273	& 268,15	& 260,396	& 7,754	& 2,891 \\
tetradecane	& 343	& 278,65	& 268,202	& 10,448	& 3,749 \\
pentadecane	& 420	& 283,05	& 275,128	& 7,922	& 2,7994 \\
hexadecane	& 512	& 291,15	& 281,901	& 9,249	& 3,177 \\
heptadecane	& 612	& 294,15	& 288,002	& 6,148	& 2,090 \\
octadecane	& 729	& 302,15	& 293,984	& 8,166	& 2,703 \\
nonadecane	& 855	& 306,15	& 299,436	& 6,714	& 2,193 \\
icosane	& 1000	& 309,85	& 304,793	& 5,057	& 1,632 \\
henicosane	& 1155	& 313,65	& 309,720	& 3,930	& 1,253 \\
docosane	& 1331	& 315,15	& 314,570	& 0,580	& 0,184 \\
tricosane	& 1518	& 322,15	& 319,066	& 3,084	& 0,957 \\
tetracosane	& 1728	& 325,15	& 323,497	& 1,653	& 0,509 \\
pentacosane	& 1950	& 327,15	& 327,630	& -0,480	& 0,147 \\
hexacosane	& 2197	& 329,55	& 331,708	& -2,158	& 0,655 \\
heptacosane	& 2457	& 332,65	& 335,533	& -2,883	& 0,867 \\
octacosane	& 2744	& 337,65	& 339,311	& -1,661	& 0,492 \\
nonacosane	& 3045	& 336,85	& 342,870	& -6,020	& 1,787 \\
triacontane	& 3375	& 338,95	& 346,388	& -7,438	& 2,195 \\
hentriacontane	& 3720	& 341,05 & 349,717	& -8,667	& 2,541 \\
dotriacontane	& 4096	& 342,15 & 353,009 & -10,859	& 3,174 \\
 \hline
average				& & & & &	4,025 
\\ \hline  
\end{tabular}
\caption{\label{tab3} Results for all 31 alkane molecules.}
\end{table}

\begin{table}
\centering
\begin{tabular}{|c|c|c|c|c|}\hline 
Molecule &	\# Aut &	$W$	& $GP$ & $MP\,(K)$ \\ \hline
1-methylnaphthalene	& 1	& 140	& 0	& -22 \\
2-methylnaphthalene	& 1	& 144	& 0	& 35 \\
1-ethylnaphthalene	& 1	& 182	& 0	& -14 \\
2-ethylnaphthalene	& 1	& 190	& 0	& -7 \\
2-6-dimethylnaphthalene	& 2	& 186	& 144	& 110 \\
2-7-dimethylnaphthalene	& 2	& 185	& 108	& 97 \\
1-7-dimethylnaphthalene	& 1 & 180	& 0	& -14 \\
1-5-dimethylnaphthalene	& 2	& 176	& 132	& 80 \\
1-2-dimethylnaphthalene	& 1	& 178	& 0	& -4 \\
2-3-5-trimethylnaphthalene	& 1	& 224	& 0	& 25 \\
anthracene	& 4	& 279	& 245	& 216 \\
1-methylanthracene	& 1	& 334	& 0	& 86 \\
2-7-dimethylanthracene	& 2	& 413	& 280	& 241 \\
9-10-dimethylanthracene	& 4	& 378	& 320	& 183 \\
phenanthrene	& 2	&271	& 175	 & 101 \\ 
3-6-dimethylphenanthrene	& 2	& 396	& 256	& 141 \\ \hline 
naphtalene	& 4	& 109	& 95	& 81 \\
1-3-7-trimethylnaphthalene	& 1	& 226	& 0	& 14 \\
2-6-dimethylanthracene	& 2	& 414	& 336	& 250 \\
4-5-methylenephenanthrene	& 2	& 300	& 165	& 116 \\ \hline 
\end{tabular}
\caption{\label{tab4} Data for the training set and the test set of PAH's.}
\end{table}

\begin{table}
\centering
\begin{tabular}{|c|c|c|c|c|c|}\hline 
Molecule & $GP$ & $MP\,(K)$ & $\widehat{MP}$	& Residual	& \% Residual \\ \hline
naphtalene &	95	& 81	& 72,686	& 8,315	& 10,265 \\
1-3-7-trimethylnaphtalene	& 0	& 14	& 10,926	& 3,074	& 21,957 \\
2-6-dimethylanthracene	& 336	& 250	& 229,360	& 20,640	& 8,256 \\
4-5-methylenephenanthrene	& 165	& 116	& 118,193	& -2,193	& 1,890 \\ \hline
average & & & &	&			10,592 \\ \hline
\end{tabular}
\caption{\label{tab5} Predicted melting points on the test set for PAH's.}
\end{table}

\begin{table}
\centering
\begin{tabular}{|c|c|c|c|}\hline 
Molecule &	\# Aut &	$GP$ &	$MP\,(K)$ \\ \hline
octane	& 2	& 64	& 216,3 \\
2-methyl-heptane	& 2	& 8	& 164,16 \\
3-methyl-heptane	& 1	& 0	& 152,6 \\
4-methyl-heptane	& 2	& 48	& 152 \\
2,2-dimethyl-hexane	& 6 & 16	& 151,97 \\
2,5-dimethyl-hexane	& 8	& 64	& 182 \\
3,3-dimethyl-hexane	& 2	& 8	& 147 \\
2-methyl-3-ethyl-pentane	& 4	& 32	& 158,2 \\
3-methyl-3-ethyl-pentane	& 6	& 48	& 182,2 \\
2,2,3-trimethyl-pentane	& 6	& 16	& 160,89 \\
2,2,4-trimethyl-pentane	& 12	& 24	& 165,8 \\
2,3,3-trimethyl-pentane	& 4	& 16	& 172,22 \\
2,3,4-trimethyl-pentane	& 8	& 48	& 163,9 \\
2,2,3,3-tetramethylbutane	& 72	& 56	& 373,8 \\ \hline  
\end{tabular}
\caption{\label{tab6} Data for octane isomers.}
\end{table}

\end{document}